\documentclass[12pt]{article}

\usepackage{latexsym}
\usepackage{amssymb}

\newtheorem{theorem}{Theorem}[section]
\newtheorem{lemma}[theorem]{Lemma}
\newtheorem{proposition}[theorem]{Proposition}
\newtheorem{corollary}[theorem]{Corollary}
\newtheorem{definition}[theorem]{Definition}
\newtheorem{example}[theorem]{Example}

\newcommand\Op{\mathop{\rm Op}}

\newcommand\nph{\varphi}

\newcommand\ccd{{\mathcal {D}}}
\newcommand\ccf{{\mathcal {F}}}

\newcommand\cce{{\mathcal {E}}}

\newcommand\ccm{{\mathcal {M}}}

\newcommand\pone{{\mathcal {P}_{1}}}
\newcommand\ptwo{{\mathcal {P}_{2}}}

\newcommand\maps{{\mathcal {M}(\pone,\ptwo)}}

\begin{document}

\title{Spectral synthesis and masa-bimodules}
\author{I. G. Todorov}
\date{ }

\maketitle

\section{Introduction and preliminaries}

The interaction between Harmonic Analysis and Operator Theory has been
fruitful. In \cite{a} Arveson related results in the Theory of Operator
Algebras to such in Spectral Synthesis. He defined synthesis for subspace
lattices and proved that certain classes of lattices are synthetic. The main
result in this note is a generalization of a result of Arveson for the case
of subspace maps. In order to describe in more detail the content of the
present work, we need to introduce some definitions and facts from \cite
{opmap} and \cite{ls}.

Let ${{\mathcal{H}}_{1}}$ and ${{\mathcal{H}}_{2}}$ be separable complex
Hilbert spaces, $\mathcal{P}_{i}$ the lattice of all (orthogonal)
projections on ${\mathcal{H}}_{i}$, $i=1,2$. Following Erdos \cite{opmap},
we let $\mathcal{M(P}_{1},\mathcal{P}_{2})$ denote the set of all maps
$\varphi :\mathcal{P}_{1}\rightarrow \mathcal{P}_{2}$ which are 0-preserving and
$\vee $-continuous (i.e. preserve arbitrary suprema). We will call such maps
\textbf{subspace maps}. It was shown in \cite{opmap} that each $\varphi \in
\mathcal{M(P}_{1},\mathcal{P}_{2})$ uniquely defines semi-lattices ${%
\mathcal{S}}_{1\varphi }\subseteq {\mathcal{P}_{1}}$ and ${\mathcal{S}}%
_{2\varphi }=\varphi ({\mathcal{P}_{1}})\subseteq {\mathcal{P}_{2}}$ such
that $\varphi $ is a bijection between ${\mathcal{S}}_{1\varphi }$ and ${%
\mathcal{S}}_{2\varphi }$ and is uniquely determined by its restriction to ${%
\mathcal{S}}_{1\varphi }$. Moreover, ${\mathcal{S}}_{1\varphi }$ is
meet-complete and contains the identity projection while ${\mathcal{S}}%
_{2\varphi }$ is join-complete and contains the zero projection. If ${%
\mathcal{S}}_{1\varphi }$ and ${\mathcal{S}}_{2\varphi }$ are commutative,
we say that $\varphi $ is a \textbf{commutative} subspace map.

Note that the set
$${\mathcal{M}}_{max}(\varphi )=\{T\in \mathcal{B}({{\mathcal{H}}_{1}},{{
\mathcal{H}}_{2}}):\varphi (P)^{\bot }TP=0\ \mbox{for each
}P\in \mathcal{P}_{1}\} $$
is also uniquely determined by $\varphi |_{\mathcal{S}_{1\varphi }}$: if $T$
satisfies $\varphi (P)^{\bot }TP=0$ for each $P\in \mathcal{S}_{1\varphi }$,
then $T\in {\mathcal{M}}_{max}(\varphi )$. Let us note that in the work of
Erdos \cite{opmap} the notation $\Op\left( \varphi \right) \ $was
used to denote the space ${\mathcal{M}}_{max}(\varphi )$. However, it was
more convenient for us to use the above notation (see below).

Given a subspace $\mathcal{U\subseteq B}({{\mathcal{H}}_{1}},{{\mathcal{H}}%
_{2}})$, we define its map $\mathop{\rm Map}\mathcal{U}:\mathcal{P}%
_{1}\rightarrow \mathcal{P}_{2}$ by
$$
\mathop{\rm Map}\mathcal{U}(P)=\overline{[\mathcal{U}(P)]}\quad (P\in
\mathcal{P}_{1})
$$
(where the symbol $\overline{[\mathcal{U}(P)]}$ stands for the projection
onto the closed subspace spanned by $\{Sx:x\in P({{\mathcal{H}}_{1}}),S\in
\mathcal{U}\}$).

If $\varphi^{\ast}=\mathop{\rm Map}{\mathcal{U}}^{\ast}$, then ${\mathcal{S}}%
_{1 \varphi}=\{P^{\perp} : P\in\varphi^* ({\mathcal{P}_{2}})\}$ \cite{opmap}.

If $\varphi_1,\dots,\varphi_n$ are maps, by $\varphi_1\wedge\dots\wedge%
\varphi_n$ we denote the map of the space ${\mathcal{M}}_{max}(\varphi_1)%
\cap\dots\cap{\mathcal{M}}_{max}(\varphi_n)$.


A subspace ${\mathcal{U}}$ is called \textbf{reflexive}, if ${\mathcal{U}}=
{\mathcal{M}}_{max}(\varphi)$ for some subspace map $\varphi$ \cite{opmap},
\cite{ls}.
If $\ccm$ is reflexive, it follows from the separability of the Hilbert spaces that
there exist (countable) families $\{E_i\}_{i=1}^{\infty}$ and $\{F_i\}_{i=1}^{\infty}$
such that ${\mathcal{M}}=\{T : F_iTE_i=0, i\in\mathbb{N}\}$.
Conversely, every set of this type is reflexive.


Given $\varphi \in {\mathcal{M}}({\mathcal{P}}_{1},{\mathcal{P}}_{2})$, the
subspace
$$
\mathcal{V}=\{T\in \mathcal{B}({{\mathcal{H}}_{1}},{{\mathcal{H}}_{2}}%
):\varphi (L)TL^{\bot }=0\ \mbox{for each
}L\in \mathcal{S}_{1\varphi }\}
$$
is clearly reflexive. We denote its map by $\varphi ^{\bot }$. Thus $%
\mathcal{V}={\mathcal{M}}_{max}(\varphi ^{\bot })$ and $\varphi ^{\bot }\in
\mathcal{M(P}_{1},\mathcal{P}_{2})$ satisfies $\varphi ^{\bot }(L^{\bot
})\leq \varphi (L)^{\bot }$ for each $L\in \mathcal{S}_{1\varphi }$.

Let us note that subspace maps are a generalization of subspace lattices.
Indeed, if ${\mathcal{L}}$ is a subspace lattice, there is a unique map
$\varphi\in\maps$ which has semi-lattices equal to ${\mathcal{L}}$ and the
restriction of $\varphi$ to ${\mathcal{L}}$ is the identity map on ${%
\mathcal{L}}$. Moreover, we have that ${\mathcal{M}}_{max}(\varphi)=%
\mathop{\rm Alg}{\mathcal{L}}$, the algebra of all operators leaving each
element of ${\mathcal{L}}$ invariant \cite{opmap}.

Now we describe the notion of synthesis for subspace maps. Recall first that
a masa-bimodule ${\mathcal{U}}\subseteq{\mathcal{B}(\mathcal{H}_{1},
\mathcal{H}_{2})}$ is a subspace of operators, such that ${\mathcal{D}_{2}}{%
\mathcal{U}}{\mathcal{D}_{1}}\subseteq{\mathcal{U}}$ for some maximal
abelian selfadjoint algebras (masas) ${\mathcal{D}_{1}}\subset{\mathcal{B}}({%
\mathcal{H}}_1)$ and ${\mathcal{D}_{2}}\subset{\mathcal{B}}({\mathcal{H}}_2)$%
. If $\varphi$ is a commutative subspace map, we consider the family
$$
\Phi=\{{\mathcal{U}}\subseteq{\mathcal{B}(\mathcal{H}_{1}, \mathcal{H}_{2})}
: {\mathcal{U}}\mbox{ is a weak*-closed  } \ccd_2,\ccd_1\mbox{-bimodule },
\mathop{\rm Map}{\mathcal{U}}
=\varphi\}.
$$
By a result of Arveson and Erdos \cite{opmap}, the family $\Phi$, along with
its maximal element, ${\mathcal{M}}_{max}(\varphi)$, contains also a minimal
element, denoted by ${\mathcal{M}}_{min}(\varphi)$. The map $\varphi$ is
called \textbf{synthetic}, if $\Phi$ contains only one element, or
equivalently, if ${\mathcal{M}}_{max}(\varphi)={\mathcal{M}}_{min}(\varphi)$%
. A reflexive masa-bimodule ${\mathcal{M}}$ is called {\bf synthetic}, if
$\mathop{\rm Map}{\mathcal{M}}$ is synthetic.

There is an alternative description of the spaces ${\mathcal{M}}%
_{max}(\varphi )$ and ${\mathcal{M}}_{min}(\varphi )$ \cite{a}, \cite{eks}.
If $(X,m)$ and $(Y,n)$ are standard measure spaces, then the sets of the
form $M\times Y\cup X\times N$, where $M\subseteq X$ and $N\subseteq Y$ are
null sets, and their subsets, are said to be \textbf{marginally null} \cite
{a}. If $\kappa \subseteq X\times Y$ is a Borel set whose complement is a
union of Borel rectangles, then an operator $T$ is said to be supported on $%
\kappa $ if $P(\beta )TP(\alpha )=0$ for each Borel rectangle $\alpha \times
\beta $ whose intersection with $\kappa $ is a marginally null set (by $%
P(\alpha )$ we denote the projection which corresponds to multiplication by
the characteristic function of the set $\alpha )$. The space ${\mathcal{M}}%
_{max}(\kappa )$ consisting of all operators, supported on $\kappa $, is a\
reflexive masa-bimodule. Each reflexive space ${\mathcal{M}}\subseteq {%
\mathcal{B}}(L^{2}(X,m),L^{2}(Y,n))$, which is a bimodule over the
multiplication masas, is of this form \cite{eks}. A set $\kappa \subseteq
X\times Y$ is called \textbf{synthetic}, if ${\mathcal{M}}_{max}(\kappa )$
is a synthetic masa-bimodule. The minimal ultraweakly closed masa-bimodule
whose map is equal to the map of ${\mathcal{M}}_{max}(\kappa )$, is denoted
by ${\mathcal{M}}_{min}(\kappa )$. To obtain more explicit description of ${%
\mathcal{M}}_{min}(\kappa )$, we need topology. If $X$ and $Y$ are compact
metric spaces with Borel measures $m$ and $n$ and if ${\mathcal{M}}\subseteq
{\mathcal{B}}(L^{2}(X,m),L^{2}(Y,n))$ is a reflexive masa-bimodule, then a
\textit{closed} set $\kappa \subseteq X\times Y$ can be chosen such that ${%
\mathcal{M}}={\mathcal{M}}_{max}(\kappa )$. Then ${\mathcal{M}}_{min}(\kappa
)$ equals to the ultaweak closure of the space of ``pseudointegral''
operators supported on $\kappa $ (see \cite{a} for details).



Recall \cite{kt} that a subspace map $\varphi$ is called an \textbf{ortho-map},
if $\varphi(L)\perp\varphi(L^{\perp})$, for each $L\in{\mathcal{S}}_{1 \chi}$%
. 
Ortho-maps are the subspace analogue of the orthocomplemented lattices. In
\cite{a} Arveson proved that each commutative Boolean lattice is synthetic. A
generalization of this fact \cite{kt} will be of importance for us:

\begin{theorem}
\label{th_norm} Commutative subspace ortho-maps are synthetic.
\end{theorem}

On the other edge of the ortho-maps are the nest maps - these are those
subspace maps which have totally ordered semi-lattices. Nest maps are also
synthetic \cite{ep}. 

Now we introduce the map analogue of the finite width lattices.

\begin{definition}
\label{d_f_w} A (commutative) subspace map $\varphi$ will be called \textbf{%
map of finite width}, if there exist nest maps $\chi_1, \chi_2,\dots,\chi_n$
such that the (respective) semi-lattices of $\chi_i$ commute with each other
and $\varphi=\chi_1\wedge\chi_2\wedge\dots\wedge\chi_n$.
\end{definition}

It is clear that if $\varphi$ is a map of finite width, then the space ${%
\mathcal{M}}_{max}(\varphi)$ is a bimodule over two finite width algebras.
As we shall see later, the converse does not hold.

\smallskip

The next Section contains the main result of the paper, namely that every
map of finite width is synthetic. In Section \ref{s_thin_bim} we investigate
the following question: Suppose that ${\mathcal{M}}$ is a masa-bimodule and
let ${\mathcal{A}}$ and ${\mathcal{B}}$ be the biggest algebras, over which $%
{\mathcal{M}}$ is a bimodule. If ${\mathcal{A}}$ and ${\mathcal{B}}$ are
synthetic, does it follow that ${\mathcal{M}}$ is synthetic as well? We
introduce a class of masa-bimodules, which we call thin masa-bimodules;
these are the reflexive masa-bimodules which are not bimodules over any
algebras bigger than the masa's. Our main example (Theorem \ref{e_non-synth})
answers the above question negatively. Next we characterize the thin
masa-bimodules in the class of normalizing masa-bimodules (for the respective definition see
the paragraph preceding Theorem \ref{p_norm_thin}). We give a criterion
for the union of two synthetic sets to be again a synthetic set and use it
to describe an example of a non-normalizing synthetic thin masa-bimodule.

\section{The main result}

\label{s_main_result}

\begin{theorem}
\label{th_synth_f_w} Maps of finite width are synthetic.
\end{theorem}

The proof of this theorem follows the ideas of Arveson \cite{a} and relies
on several lemmas.

\begin{lemma}
\label{l_supports} Let $\{{\mathcal{M}}_{i}\}_{i=1}^{n}$ (${\mathcal{M}}%
_{i}\subseteq {\mathcal{B}(\mathcal{H}_{1},\mathcal{H}_{2})}$) be a (finite)
family of reflexive masa-bimodules, let $\varphi _{i}=%
\mathop{\rm
Map}{\mathcal{M}}_{i}$, $i=1,\dots ,n$, and suppose that the respective
semi-lattices of $\varphi _{i}$, $i=1,\dots ,n$, commute with each other.
Then there are compact metric spaces $X$ and $Y$, regular Borel measures $m$
and $n$ on $X$ and $Y$ respectively and closed sets $\kappa _{i}\subseteq
X\times Y$ such that ${\mathcal{H}}_{1}$ and ${\mathcal{H}}_{2}$ are
unitarily equivalent to $L^{2}(X,m)$ and $L^{2}(Y,n)$, the masas are
unitarily equivalent to the multiplication algebras of $L^{\infty }(X,m)$
and $L^{\infty }(Y,n)$ and ${\mathcal{M}}_{i}$ is unitarily equivalent to ${%
\mathcal{M}}_{max}(\kappa _{i})$, for each $i=1,\dots ,n$. Moreover, if ${%
\mathcal{M}}=\cap _{i=1}^{n}{\mathcal{M}}_{i}$, then ${\mathcal{M}}={%
\mathcal{M}}_{max}(\cap _{i=1}^{n}\kappa _{i})$.
\end{lemma}

\noindent The lemma follows immediately from the work in Theorem 4.6 of \cite{eks} and we
omit the proof.

\begin{lemma}
\label{l_rel} If $\nph$ and $\psi$ are commutative maps and $
{\mathcal{M}}_{max}(\varphi )$ $\subseteq $ ${\mathcal{M}}_{max}(\psi )$,
then ${\mathcal{M}}_{min}(\varphi )$ $\subseteq $ ${\mathcal{M}}_{min}(\psi
) $.
\end{lemma}

\noindent \textit{Proof. } From Lemma \ref{l_supports} it follows that we
can choose representations ${\mathcal{H}}_1=L^2(X,m)$ and ${\mathcal{H}}%
_2=L^2(Y,n)$ and closed sets $\kappa_1$ and $\kappa_2$ with $%
\kappa_1\subseteq\kappa_2$ such that ${\mathcal{M}}_{max}(\varphi_1)={%
\mathcal{M}}_{max}(\kappa_1)$ and ${\mathcal{M}}_{max}(\varphi_2)={\mathcal{M%
}}_{max}(\kappa_2)$. From \cite{a} we have that ${\mathcal{M}}%
_{min}(\kappa_1)$ is the ultraweak closure of the space of pseudointegral
operators with measures supported on $\kappa_1$, and similarly for ${%
\mathcal{M}}_{min}(\kappa_2)$. The conclusion follows since $%
\kappa_1\subseteq\kappa_2$. $\diamondsuit $\hfill \bigskip


\begin{lemma}
\label{l_synth_nest_map} Let $\varphi$ be a synthetic commutative subspace map and
$\chi$ be a nest map with semi-lattices $\Sigma_1$ and $\Sigma_2$ which
commute with the semi-lattices of the map $\varphi$. Then, for each $%
P\in\Sigma_1$, we have
$$
\chi(P){\mathcal{M}}_{max}(\varphi)P^{\perp}\subseteq{\mathcal{M}}%
_{min}(\varphi\wedge\chi).
$$
\end{lemma}

\noindent \textit{Proof. } Choose
families $\cce=\{E_i\}_{i=1}^{\infty}\subseteq \Sigma_1$ and
$\ccf=\{F_i\}_{i=1}^{\infty}\subseteq\Sigma_2$ such that
$P\in \{E_i\}_{i=1}^{\infty}$,
$\chi(E_i) = F_i$, $\cce$ is strongly dense in $\Sigma_1$ and
$\ccf$ is strongly dense in $\Sigma_2$.
Next, acording to Theorem  4.6 of \cite{eks} and Lemma \ref{l_supports} choose
representations ${\mathcal{H}}_1=L^2(X,m)$,
${\mathcal{H}}_2=L^2(Y,n)$,
in such a way that there exist closed and open sets $\alpha_i\subseteq X$
and $\beta_i\subseteq Y$ such that $E_i = P(\alpha_i)$ and
$F_i = P(\beta_i)$ and the support $\kappa$ of $\ccm_{max}(\nph)$ can be taken to be
a closed set.
Then the support  $\lambda\subseteq X\times Y$
is (marginally equivalent to) the complement of the set
$\cup_{i=1}^{\infty} \alpha_i\times\beta_i^c$.
Let $\alpha=\alpha_{i_0}$ and
$\beta=\beta_{i_0}$ where $i_0$ is the index with the property
$P=P(\alpha_{i_0})$.
It is straightforward to verify that $\alpha^c\times\beta\subset\lambda$.

Let $T\in{\mathcal{M}}_{max}(\varphi)$. Since $\varphi$ is synthetic, it
follows that $T$ is a limit of a net $\{T_{\gamma}\}_{\gamma\in\Gamma}$ of
pseudointegral operators supported on $\kappa$. Then the operators
$S_{\gamma}=\chi(P)T_{\gamma}P^{\perp}$, $\gamma\in\Gamma$, are
pseudointegral and the previous paragraph shows that $S_{\gamma}$
is supported on $\kappa\cap\lambda$, for each $\gamma\in\Gamma$. Obviously $
S_{\gamma}\longrightarrow \chi(P)TP^{\perp}$. This shows that $
\chi(P)TP^{\perp}\in{\mathcal{M}}_{min}(\kappa\cap\lambda)$. $\diamondsuit $
\hfill \bigskip


\begin{lemma}
\label{l_two_synth_maps} Let $\varphi$ be a commutative subspace map and $\chi$ be a
nest map such that the respective semi-lattices of $\varphi$ and $\chi$
commute. If the maps $\varphi$ and $\varphi\wedge\chi\wedge\chi^{\perp}$ are
synthetic, then the map $\varphi\wedge\chi$ is synthetic.
\end{lemma}

\noindent \textit{Proof. } Let $T\in{\mathcal{M}}_{max}(\varphi\wedge\chi)$.
We will show that $T\in{\mathcal{M}}_{min}(\varphi\wedge\chi)$. Let $%
\Sigma_1 $ and $\Sigma_2$ be the semi-lattices of $\chi$. Since $\chi$ is a
nest map, $\Sigma_1$ and $\Sigma_2$ are totally ordered. Let $%
\{E_i\}\subseteq\Sigma_1$ be a countable family dense in $\Sigma_1$ in the
strong operator topology. If $E$ is a projection, we set $E^{+1}=E$ and $%
E^{-1}=E^{\perp}$. For each $s=(s_1,\dots,s_n)$, where $n\in\mathbb{N}$,
we let
$$
E^s=E_1^{s_1}E_2^{s_2}\dots E_n^{s_n}.
$$
It is clear that the projections $E^s$ are just the atoms of the Boolean
lattice generated by the family $\{E_i\}_{i=1}^n$. We set also
$$
\chi(E^s)= \chi(E_1)^{s_1}\chi(E_2)^{s_2}\dots \chi(E_n)^{s_n}.
$$
We let
$$
R_n=\sum_s \chi(E^s)TE^s, \ \ S_n=\sum_{s\neq t}\chi(E^s)TE^t.
$$
Obviously $T=R_n+S_n$, for each $n\in\mathbb{N}$. We will show first that $%
S_n\in{\mathcal{M}}_{min}(\varphi\wedge\chi)$. It is enough to see that $%
\chi(E^s)TE^t\in{\mathcal{M}}_{min}(\varphi\wedge\chi)$, if $s\neq t$. So
let $s\neq t$ and $i$, $1\leq i\leq n$, be such that $s_i\neq t_i$. If $%
s_i=+1$ and $t_i=-1$, then there are projections $F$ and $G$ in the
respective masas such that $\chi(E^s)TE^t=G\chi(E_i)^{\perp}TE_iF$. The last
expression is equal to zero, since $T\in{\mathcal{M}}_{max}(\varphi\wedge%
\chi)\subseteq{\mathcal{M}}_{max}(\chi)$. If $s_i=+1$ and $t_i=-1$, then
there are projections $F$ and $G$ in the respective masas such that $%
\chi(E^s)TE^t=G\chi(E_i)TE_i^{\perp}F$. From Lemma \ref{l_synth_nest_map} it
follows that $\chi(E_i)TE_i^{\perp}\in{\mathcal{M}}_{min}(\varphi\wedge\chi)$
and thus $\chi(E^s)TE^t\in{\mathcal{M}}_{min}(\varphi\wedge\chi)$ and so we
proved that $S_n\in{\mathcal{M}}_{min}(\varphi\wedge\chi)$.

In order to deal with $R_n$, we notice first that for each $i=1,\dots,n$ we
have that $R_nE_i=\chi(E_i)R_n$. Since the sequence $\{R_n\}_{n=1}^{\infty}$
is $\|\cdot\|$-bounded, there is a subsequence $\{R_{n^{\prime}}\}$,
ultraweakly converging to an operator $R$. It is easy to see that $R\in{%
\mathcal{M}}(\chi\wedge\chi^{\perp})$. Indeed, for each $n^{\prime}$ and $%
m^{\prime}$ with $m^{\prime}>n^{\prime}$, we have that $\chi(E_i)R_{m^{%
\prime}}=R_{m^{\prime}}E_i$, $i=1,\dots,n^{\prime}$. Thus $\chi(E_i)R=RE_i$
for each $i=1,\dots,n^{\prime}$. Since this holds for each $n^{\prime}$, we
conclude that $R\in{\mathcal{M}}_{max}(\chi\wedge\chi^{\perp})$.

If $S=T-R$, we have that $S_{n^{\prime }}\longrightarrow S$ in the ultraweak
topology. We showed that $S_{n^{\prime }}\in {\mathcal{M}}_{min}(\varphi
\wedge \chi )$. Since the last space is ultarweakly closed, we have that $%
S\in {\mathcal{M}}_{min}(\varphi \wedge \chi )$. On the other hand, we have
that
$$
{\mathcal{M}}_{max}(\varphi \wedge \chi \wedge \chi ^{\perp })\subset
{\mathcal{M}}_{max}(\varphi \wedge \chi ).
$$
So, according to Lemma \ref{l_rel}, we have that
$$
R\in {\mathcal{M}}_{max}(\varphi \wedge \chi \wedge \chi ^{\perp })={%
\mathcal{M}}_{min}(\varphi \wedge \chi \wedge \chi ^{\perp })\subset {%
\mathcal{M}}_{min}(\varphi \wedge \chi ).
$$
Thus $T=R+S\in {\mathcal{M}}_{min}(\varphi \wedge \chi )$. $\diamondsuit $%
\hfill \bigskip

\bigskip

\noindent \textit{Proof of Theorem \ref{th_synth_f_w}. } It suffices to show
the following:

\noindent \textbf{Claim.} If $\chi_1,\dots,\chi_n$ are nest maps, $%
\psi=\chi_1\wedge \chi_1^{\perp}\wedge\dots\wedge\chi_n\wedge\chi_n^{\perp}$%
, $\varphi$ is a finite width map such that the respective semi-lattices of $%
\chi_i$ and $\varphi$ commute with each other, then the map $%
\psi\wedge\varphi$ is synthetic.

Indeed, suppose that the claim holds and let $\varphi$ be a finite width
map. We will show that $\varphi$ is synthetic using induction on its width.
If $\varphi$ is a nest map, then it is of course synthetic. Now let $\varphi$
have width $n+1$, $\varphi_0$ be a map of width $n$ and $\chi$ a nest map
such that $\varphi=\varphi_0\wedge\chi$. From the Claim we have that the map
$\varphi_0\wedge\chi\wedge\chi^{\perp}$ is synthetic. Now from Lemma \ref
{l_two_synth_maps} it follows that $\varphi_0\wedge\chi=\varphi$ is
synthetic.

To prove the Claim we use again induction on the width of $\varphi $. Let $%
\varphi $ be a nest map. From \cite{kt} we have that $\psi $ and $\psi
\wedge \varphi \wedge \varphi ^{\perp }$ are synthetic maps. From Lemma \ref
{l_two_synth_maps} we conclude that $\psi \wedge \varphi $ is synthetic. Let
$\varphi $ have width $n+1$. Then $\varphi =\varphi _{0}\wedge \chi $ where $%
\varphi _{0}$ has width $n$ and $\chi $ is a nest map. Thus $\varphi
_{0}\wedge \psi $ and $\varphi _{0}\wedge \psi \wedge \chi \wedge \chi
^{\perp }$ are synthetic and so again from Lemma \ref{l_two_synth_maps} we
obtain that $\varphi _{0}\wedge \psi \wedge \chi =\varphi \wedge \psi $ is
synthetic. $\diamondsuit $\hfill \bigskip

\noindent\textbf{Remarks (i)} We do not know whether the result on the
synthesis of finite width maps can be inferred directly from the Arveson's
result for finite width lattices by using a suitable two by two matrix
algebra.

\smallskip

\noindent \textbf{(ii)} It is known that the support of a nest algebra bimodule is of
the form $\left\{ (x,y):f(x)\leq g(y)\right\} $, for some Borel functions $f$
and $g$ (see \cite{t} for a detailed proof). Thus Theorem \ref{th_synth_f_w}
can be rephrased by saying that the set of solutions of a system of
inequalities of the form
$$
f_{1}(x)\leq g_{1}(y)
$$
$$
f_{2}(x)\leq g_{2}(y)
$$
$$
.................
$$
$$
f_{n}(x)\leq g_{n}(y),
$$
where $f_{1},f_{2},\dots ,f_{n}$ and $g_{1},g_{2},\dots ,g_{n}$ are Borel
functions, is a synthetic set.

\smallskip

\noindent \textbf{(iii) }Theorem \ref{th_synth_f_w} is a generalization of the result
on the synthesis of orthomaps (Theorem \ref{th_norm}) \cite{kt}, \cite{shulman}.
Indeed, a set of the form
$$\left\{ (x,y):f(x)=g(y)\right\}$$
is clearly the set of solutions of a system of two inequalities
(see the paragraph before Theorem \ref{p_norm_thin}).

\section{Thin masa-bimodules}

\label{s_thin_bim}

If a reflexive subspace ${\mathcal{M}}\subseteq {\mathcal{B}(\mathcal{H}_{1},%
\mathcal{H}_{2})}$ is a bimodule over two masas ${\mathcal{D}_{1}}\subseteq {%
\mathcal{B}}({\mathcal{H}}_{1})$ and ${\mathcal{D}_{2}}\subseteq {\mathcal{B}%
}({\mathcal{H}}_{2})$, then ${\mathcal{M}}$ is automatically a bimodule over
two, in general bigger, CSL algebras ${\mathcal{A}}$ and ${\mathcal{B}}$.
Indeed, if $\varphi =\mathop{\rm Map}{\mathcal{M}}$ and ${\mathcal{S}}_{1}$
and ${\mathcal{S}}_{2}$ are the semi-lattices of $\varphi $, then it is
enough to put ${\mathcal{A}}=\mathop{\rm Alg}{\mathcal{S}}_{1}$ and ${%
\mathcal{B}}=\mathop{\rm Alg}{\mathcal{S}}_{2}$. It is also easy to see that
the algebras ${\mathcal{A}}$ and ${\mathcal{B}}$ are the biggest algebras
over which the space ${\mathcal{M}}$ is bimodule. Of course, crucial here is
the fact that ${\mathcal{M}}$ is reflexive. We are interested in the
following question: If the algebras ${\mathcal{A}}$ and ${\mathcal{B}}$ are
synthetic, does it follow that the bimodule ${\mathcal{M}}$ is synthetic as
well? We show that the answer to this question is negative. More generally,
we investigate the class of masa-bimodules which are not bimodules over any
bigger algebras.

\begin{definition}
\label{d_thin} A reflexive masa-bimodule ${\mathcal{M}}$ will be called
\textbf{thin}, if the biggest algebras over which ${\mathcal{M}}$ is a
bimodule, are the masas.
\end{definition}

It is clear that a masa is a thin masa-bimodule. A simple class of thin masa-bimodules
consists of the bimodules of the form $U\ccd V$, where $\ccd$ is a masa and $U$ and
$V$ are unitary operators (the last space is obviously a bimodule over the masas $\ccd_1 = U\ccd U^*$
and $\ccd_2 = V^*\ccd V$). It is not true, however, that every thin masa-bimodule is of this
form (see Theorem \ref{e_non-synth} and Example \ref{e_nn_sy_th}).

We first give a criterion which we will use in determining whether a (reflexive)
masa-bimodule is thin.
It is an immediate consequence of a rather well-known condition that a family of
projections generates the multiplication masa (as a W*-algebra).
The proof is omited.

\begin{proposition}
\label{p_cr} Let $(Z,\mu )$ be a standard measure space, ${\mathcal{H}}=L^{2}(Z,\mu )$
and ${\mathcal{D}}$ be the multiplication masa. Let ${\mathcal{S}}\subseteq
{\mathcal{D}}$ be a family of projections and ${\mathcal{C}}=
\mathop{\rm Alg}{\mathcal{S}}$. Then ${\mathcal{C}}={\mathcal{D}}$ if and only if there
is a countable family $\{\beta _{i}\}_{i=1}^{\infty }$ of Borel subsets of $Z$
such that $\{P(\beta _{i})\}_{i=1}^{\infty }\subseteq {\mathcal{S}}$ and
which separates almost all points of $Z$ (in the sense that there is a null
set $A\subseteq Z$ such that for each $x,y\not\in A$, there is an $i$ such that
$x\in \beta _{i}$ and $y\not\in \beta _{i}$).
\end{proposition}

Before proceeding, we state two Propositions which will be useful for us
throughout this section. They relate the work of Erdos \cite{opmap} to that
of Erdos, Katavolos and Shulman \cite{eks}. \ The proof of the first
Proposition is omited. Let $(X,m)$ and $(Y,n)$ be standard measure spaces,
${\mathcal{H}}_{1}=L^{2}(X,m)$ and ${\mathcal{H}}_{2}=L^{2}(Y,n)$.

\begin{proposition}
\label{p_map_supp} Let $\kappa\subseteq X\times Y$ be a Borel set, ${%
\mathcal{U}}={\mathcal{M}}_{max}(\kappa)$ and $\varphi=%
\mathop{\rm
Map}{\mathcal{U}}$. Let $\alpha\subseteq X$ and $\beta\subseteq Y$. Then $%
\varphi(P(\alpha))=P(\beta)^{\perp}$ if and only if

(a) $(\alpha\times\beta)\cap\kappa\simeq\emptyset$ and

(b) if $(\alpha\times\beta^{\prime})\cap\kappa\simeq\emptyset$ for some
Borel set $\beta^{\prime}\subseteq Y$, then $\beta^{\prime}\subseteq\beta$
modulo null set.
\end{proposition}

\noindent

\begin{proposition}
\label{p_norm_map_supp} Let $f:Y\longrightarrow X$ be a Borel function, $%
\kappa =\{(f(y),y):y\in Y\}$, ${\mathcal{U}}={\mathcal{M}}_{max}(\kappa )$
and $\varphi =\mathop{\rm Map}{\mathcal{U}}$. Then $\varphi (P(\alpha
))=P(f^{-1}(\alpha ))$ if and only if $f$ preserves null sets (in the sense
that $f^{-1}(\alpha )$ is null whenever $\alpha $ is null).
\end{proposition}

\noindent \textit{Proof. } If $\varphi$ is given as in the Lemma and $%
\alpha\subseteq X$ is a null set, then $P(\alpha)=0$ and thus $%
\varphi(P(\alpha))=0$. This means that $P(f^{-1}(\alpha))=0$ and so $%
f^{-1}(\alpha)$ is a null set.

Conversely, suppose that $f$ preserves null sets. Fix a Borel set $\alpha
\subseteq X$. It is clear that $(\alpha \times f^{-1}(\alpha )^{c})\cap
\kappa =\emptyset $. Let $\beta \subseteq Y$ be such that $(\alpha \times
\beta )\cap \kappa \simeq \emptyset $. Put $\gamma =\beta \cap f^{-1}(\alpha
)$. We have that
$$
\{(f(y),y):y\in \gamma \}\subseteq (\alpha \times \beta )\cap \kappa
$$
and thus
$$
\{(f(y),y):y\in \gamma \}\subseteq M\times Y\cup X\times N,
$$
for some null sets $M\subset X$ and $N\subset Y$. This means that, if $y\in
\gamma $, then either $y\in N$ or $f(y)\in M$. Thus $f(\gamma \cap
N^{c})\subseteq M$, or, $\gamma \cap N^{c}\subseteq f^{-1}(M)$. Since $f$
preserves null sets, we have that $\gamma \cap N^{c}$ is null, and so $%
\gamma $ is null. This means that $\beta \subseteq f^{-1}(\alpha )^{c}$
modulo null set and from Proposition \ref{p_map_supp} we have that $%
P(f^{-1}(\alpha ))=\varphi (P(\alpha ))$. $\diamondsuit $\hfill \bigskip

The main result in this section is the following.

\begin{theorem}
\label{e_non-synth} There exists a non-synthetic thin masa-bimodule.
\end{theorem}

\noindent \textit{Proof. } Let $\mathbb{T}=[0,1)$ be the group of the circle
and let $E\subset\mathbb{T}$ be a closed set such that

(a) $E$ fails spectral synthesis;

(b) $\mathop{\rm diam} E <\frac{1}{3}$;

(c) if $a=\inf E$ and $b=\sup E$, then $\frac{1}{2}$ is the centre of the
interval $[a,b]$.

\noindent In order to establish the existence of such an $E$ we recall that
each Cantor set of $\mathbb{T}$ contains a set which fails spectral
synthesis (see \cite{katz}, p. 231). We choose a Cantor set $C$ with $0=\inf
C$ and $\mathop{\rm diam}C<\frac{1}{3}$ and a subset $E_{0}\subseteq C$
which fails spectral synthesis. Then we translate $E_{0}$ to get the
required set $E$. (We note that the class of sets of spectral synthesis is
invariant under translations.)

Let $s=\frac{b-a}{2}$ and $\lambda=\{(x,y)\in\mathbb{T}\times\mathbb{T} :
y-x\in E\}$. It follows from \cite{f} that the space ${\mathcal{M}}={%
\mathcal{M}}_{max}(\lambda) \subseteq{\mathcal{B}}(L^2(\mathbb{T}))$ is a
non-synthetic masa-bimodule. In order to finish the proof we need only show
that ${\mathcal{M}}$ is thin. Let $\varphi=\mathop{\rm Map} {\mathcal{M}}$
and ${\mathcal{S}}_1$ and ${\mathcal{S}}_2$ be the semi-lattices of $\varphi$%
. We will show only that $\mathop{\rm Alg}{\mathcal{S}}_2$ is a masa, the
fact that $\mathop{\rm Alg}{\mathcal{S}}_1$ is also a masa obviously follows
by symmetry. We will apply Proposition \ref{p_cr}. Let $\{r_i\}_{i\in\mathbb{%
N}}$ be the set of rationals in $\mathbb{T}$. We fix $i\in\mathbb{N}$.

\smallskip

\noindent 1) If $0<r_i<a-2s$, we let $\alpha_i^1=[0,a+r_i]$, $%
\alpha_i^2=[b+r_i,1]$, $\beta_i^1=[r_i,a]$ and $\beta_i^2=[0,r_i]\cup [b,1]$%
. Then according to Proposition \ref{p_map_supp}, we have that $%
\varphi(P(\alpha_i^j))^{\perp}=P(\beta_i^j)$, $j=1,2$.

\smallskip

\noindent 2) If $a-2s<r_i<a$, we let $\alpha_i^1=[0,a+r_i]$, $%
\alpha_i^2=[b+r_i,1]$, $\beta_i^1=[r_i,a]$, $\gamma_i^2=[0,r_i]\cup [b,1]$
and $\delta_i=[0,r_i]\cup [b,1]\cup [a, r_i+2s]$. By Proposition \ref
{p_map_supp}, we have that $P(\beta_i^1)=\varphi(P(\alpha_i^1))^{\perp}$. If
$\beta_i^2$ is such that $P(\beta_i^2)=\varphi(P(\alpha_i^2))^{\perp}$, then
it is easy to see that $\gamma_i^2\subseteq\beta_i^2\subseteq \delta_i^2$,
modulo null sets.

\smallskip

\noindent 3) If $a<r_{i}<b$, we let $\alpha _{i}^{1}=[r_{i}-a,b]$, $\alpha
_{i}^{2}=[a,a+r_{i}]$, $\beta _{i}^{1}=[2s,r_{i}]$ and $\beta
_{i}^{2}=[r_{i},1-2s]$. Again by Proposition \ref{p_map_supp}, we have that $%
\varphi (P(\alpha _{i}^{j}))^{\perp }=P(\beta _{i}^{j})$, $j=1,2$.

\smallskip

\noindent 4) If $b<r_{i}<b+2s$, we let $\alpha _{i}^{1}=[0,r_{i}-b]$, $%
\alpha _{i}^{2}=[r_{i}-a,1]$, $\gamma _{i}^{1}=[0,a]\cup \lbrack r_{i},1]$, $%
\delta _{i}^{1}=[0,a]\cup \lbrack r_{i}-2s,b]\cup \lbrack r_{i},1]$, $\beta
_{i}^{2}=[b,r_{i}]$. As above we have that $P(\beta _{i}^{1})=\varphi
(P(\alpha _{i}^{1}))^{\perp }$. If $\beta _{i}^{2}$ is such that $P(\beta
_{i}^{2})=\varphi (P(\alpha _{i}^{2}))^{\perp }$, then it is easy to see
that $\gamma _{i}^{2}\subseteq \beta _{i}^{2}\subseteq \delta _{i}^{2}$,
modulo null sets.

\smallskip

\noindent 5) If $b+2s<r_{i}<1$, we let $\alpha _{i}^{1}=[0,r_{i}-b]$, $%
\alpha _{i}^{2}=[r_{i}-a,1]$, $\beta _{i}^{1}=[0,a]\cup \lbrack r_{i},1]$
and $\beta _{i}^{2}=[b,r_{i}]$. As above we have $\varphi (P(\alpha
_{i}^{j}))^{\perp }=P(\beta _{i}^{j})$, $j=1,2$.

\smallskip

We have that $\{P(\beta _{i}^{jc}):i\in \mathbb{N},j=1,2\}\subseteq {%
\mathcal{S}}_{2}$, since the projections $P(\beta _{i}^{j})^{\perp }$ are in
the image of $\varphi $. In order to complete the proof, it will suffice to
show that the family $\{\beta _{i}^{j}:i\in \mathbb{N},j=1,2,3\}$ separates
almost all points of $\mathbb{T}$, since then, by symmetry, the family $%
\{\beta _{i}^{jc}:i\in \mathbb{N},j=1,2,3\}$ will also separate almost all
points of $\mathbb{T}$. So we will show that $\{\beta _{i}^{j}:i\in \mathbb{N%
},j=1,2\}$ separates the points of $\mathbb{T}\setminus
\{0,a-2s,a,b,b+2s,1\} $. Let $x\in \mathbb{T}\setminus \{0,a-2s,a,b,b+2s,1\}$%
. We consider the following cases:

\smallskip

\noindent 1) $0<x<a-2s$. Then
\begin{eqnarray*}
\{x\} &\subseteq &\cap \{\beta _{i}^{j}:x\in \beta _{i}^{j},i\in \mathbb{N}%
,j=1,2\} \\
&\subseteq &\cap \{\beta _{i}^{j}:x\in \beta _{i}^{j},0<r_{i}<a-2s,j=1,2\} \\
&=&\cap \{\beta _{i}^{1}\cap \beta _{k}^{2}:x\in \beta _{i}^{1}x\in \beta
_{k}^{2},0<r_{k}<r_{i}<a-2s\} \\
&=&\cap \{[r_{k},r_{i}]:x\in \lbrack r_{k},r_{i}]\}=\{x\}.
\end{eqnarray*}

\smallskip

\noindent 2) $a-2s<x<a$. Then
\begin{eqnarray*}
\{x\} &\subseteq &\cap \{\beta _{i}^{j}:x\in \beta _{i}^{j},i\in \mathbb{N}%
,j=1,2\} \\
&\subseteq &\cap \{\beta _{i}^{j}:x\in \beta _{i}^{j},a-2s<r_{i}<a,j=1,2\} \\
&=&\cap \{\beta _{i}^{1}\cap \beta _{k}^{2}:x\in \beta _{i}^{1}x\in \beta
_{k}^{2},a-2s<r_{i}<r_{k}<a\}.
\end{eqnarray*}
Since $\gamma _{k}^{2}\subseteq \beta _{k}^{2}\subseteq \delta _{k}^{2}$, we
have that $\beta _{i}^{1}\cap \gamma _{k}^{2}\subseteq \beta _{i}^{1}\cap
\beta _{k}^{2}\subseteq \beta _{i}^{1}\cap \delta _{k}^{2}$. If $r_{i}<r_{k}$%
, then $\beta _{i}^{1}\cap \gamma _{k}^{2}=\beta _{i}^{1}\cap \delta
_{k}^{2}=[r_{i},r_{k}]\cup \{a\}$. This implies that $\beta _{i}^{1}\cap
\beta _{k}^{2}=(r_{i},r_{k})\cup \{a\}$ and it follows that $\cap \{\beta
_{i}^{j}:x\in \beta _{i}^{j},i\in \mathbb{N},j=1,2\}\setminus \{a\}=\{x\}$.

\smallskip

\noindent 3) $a<x<b$. Then
\begin{eqnarray*}
\{x\} & \subseteq &\cap\{\beta_i^j : x\in\beta_i^j, i\in\mathbb{N}, j=1,2\}
\subseteq \cap\{\beta_i^j : x\in\beta_i^j, a<r_i<b, j=1,2\} \\
& = & \cap\{\beta_i^1\cap\beta_k^2 : x\in\beta_i^1 x\in\beta_k^2,
a<r_i<r_k<b\} \\
& = & \cap\{[r_i,r_k] : x\in [r_i,r_k]\}=\{x\}.
\end{eqnarray*}

\smallskip

\noindent 4) $b<x<b+2s$. This case is similar to 2)

\smallskip

\noindent 5) $b+2s<x<1$. This case is similar to 1) $\diamondsuit $\hfill
\bigskip

\noindent {\bf Remarks}
{\bf (i)} It is clear that the only thin (reflexive) masa-bimodule that is also an algebra, is
the masa. Thus, one might imagine that thin (reflexive) masa-bimodules are the module
analogue of the masa's. Theorem \ref{e_non-synth} shows that this is not the case.
The proper module analogue of the masa are the {\it normalizing} thin masa-bimodules,
which are characterized in Theorem \ref{p_norm_thin}.

\smallskip

\noindent {\bf (ii)} Every thin CSL algebra is synthetic since it must be equal to the masa. Thus
Theorem \ref{e_non-synth} contains another example where algebra results do not
carry over to bimodules.

Theorem \ref{e_non-synth} implies the following.

\begin{corollary}
\label{c_f_m} There exists a reflexive masa-bimodule ${\mathcal{M}}$ such
that the maximal algebras over which it is a bimodule are of finite width,
but ${\mathcal{M}}$ is not a finite width bimodule.
\end{corollary}

\noindent \textit{Proof. } Let ${\mathcal{M}}$ be the masa-bimodule from the
proof of Theorem \ref{e_non-synth}. Then the maximal algebras over which ${%
\mathcal{M}}$ is a bimodule are masas, and so have width two.
On the other hand, $\ccm$ is not a finite width bimodule itself, because if it
were, then it would have been synthetic (Theorem \ref{th_synth_f_w}).
$\diamondsuit $\hfill \bigskip

We proceed by describing a class of examples of thin masa-bimodules. Let us
recall that an ultraweakly closed masa-bimodule ${\mathcal{U}}$ is called
\textbf{normalizing} \cite{kt}, if $\mathop{\rm Map}{\mathcal{U}}$ is an ortho-map.
(These are precisely the reflexive masa-bimodules which are also {\it ternary rings of
operators} in the sense of \cite{zettl}, that is, closed under the triple product
$(T,S,R)\longrightarrow TS^*R$. Here we have retained the terminology from
\cite{kt}.)
The ultraweakly closed normalizing masa-bimodules were characterized in \cite{kt}
as the reflexive masa-bimodules whose support have the form
$\{(x,y):f(x)=g(y)\}$, where $f:X\longrightarrow \mathbb{R}$ and
$g:Y\longrightarrow \mathbb{R}$ are measurable functions.

\begin{theorem}
\label{p_norm_thin} Let $(X,m)$ and $(Y,n)$ be standard measure spaces,
${\mathcal{H}}_{1}=L^{2}(X,m)$, ${\mathcal{H}}_{2}=L^{2}(Y,n)$, ${\mathcal{D}%
_{1}}$ and ${\mathcal{D}_{2}}$ be the multiplication masas on ${\mathcal{H}}%
_{1}$ and ${\mathcal{H}}_{2}$. Then an ultraweakly closed normalizing ${%
\mathcal{D}_{2}},{\mathcal{D}_{1}}$-bimodule ${\mathcal{U}}$ is thin
if and only if there exists an almost one-to-one and onto Borel
function $f:X\longrightarrow Y$ which preserves null sets, such that ${%
\mathcal{U}}={\mathcal{M}}_{max}(\kappa )$, where $\kappa =\{(x,f(x)):x\in
X\}$ is the graph of $f$.
\end{theorem}

\noindent \textit{Proof. } Suppose that a normalizing masa-bimodule ${%
\mathcal{U}}$ is thin. Since ${\mathcal{U}}$ is a bimodule over the algebras
$({\mathcal{U}}{\mathcal{U}}^{\ast })^{\prime \prime }$ and $({\mathcal{U}}%
^{\ast }{\mathcal{U}})^{\prime \prime }$, it follows that these algebras are
abelian. By Corollary 3.3 of \cite{kt} we have that there exist Borel sets $%
X_{0}\subseteq X$ and $Y_{0}\subseteq Y$ and Borel functions $%
f:X_{0}\longrightarrow Y$ and $g:Y_{0}\longrightarrow X$ which preserve null
sets, such that both of the sets
$$
\{(g(y),y):y\in Y_{0}\}
$$
and
$$
\{(x,f(x)):x\in X_{0}\}
$$
are marginally equivalent to the support of ${\mathcal{U}}$. It follows that
$f$ and $g$ are (almost) one-to-one functions. On the other hand, ${\mathcal{%
U}}$ is easily seen to be a bimodule over the algebras $Q{\mathcal{B}}({%
\mathcal{H}}_{2})Q$ and $P{\mathcal{B}}({\mathcal{H}}_{1})P$, where $P$ and $%
Q$ are the projections corresponding to the complements of $X_{0}$ and $%
Y_{0} $ respectively. Since these algebras must be abelian, it follows that $%
P=Q=0$, so that $X_{0}$ and $Y_{0}$ have full measure. Thus the support of ${%
\mathcal{U}}$ is the graph of an almost one-to-one and onto function which
preserves null sets.

Conversely, suppose that $f:X\longrightarrow Y$ is an almost one-to-one and
onto Borel function which preserves null sets and let ${\mathcal{U}}={%
\mathcal{M}}_{max}(\kappa )$, where $\kappa =\{(x,f(x)):x\in X\}$. By
Proposition \ref{p_norm_map_supp} we have that $\varphi (P(\alpha
))=P(f(\alpha ))$, for each Borel set $\alpha \subseteq X$. Thus, for each
Borel set $\alpha \subseteq X$, the projection $P(f(\alpha ))$ is an element
of the right semi-lattice of $\varphi $. Since $X$ is standard, there is a
family $\{\alpha _{i}\}_{i=1}^{\infty }$ which separates the points of $X$.
Since $f$ is one-to-one and onto, the family $\{f(\alpha
_{i})\}_{i=1}^{\infty }$ separates the points of $Y$. From Proposition \ref
{p_cr} it follows that the maximal algebra over which ${\mathcal{U}}$ is
right module is the masa. By symmetry we get that ${\mathcal{U}}$ is a thin
masa-bimodule. $\diamondsuit $\hfill \bigskip

\bigskip

Before giving the last example of this section, we want to state a simple
criterion, in terms of masa-bimodules, for the union of two synthetic sets
to be again a synthetic set. We state also two Corollaries of this
criterion. Let $X$ and $Y$ be compact metric spaces and $m$ and $n$
regular Borel measures on $X$ and $Y$ respectively.
Let ${\mathcal{H}}_1=L^2(X,m)$ and ${\mathcal{H}}_2=L^2(Y,n)$.

\begin{proposition}
\label{p_u_s} Let $\kappa_1,\kappa_2 \subseteq X\times Y$ be closed
synthetic sets. Then the set $\kappa=\kappa_1\cup \kappa_2$ is synthetic if
and only if ${\mathcal{M}}_{max}(\kappa_1)+{\mathcal{M}}_{max}(\kappa_2)$ is
dense in ${\mathcal{M}}_{max}(\kappa)$ in the ultarweak topology.
\end{proposition}

\noindent \textit{Proof. } First we show that
\begin{eqnarray*}
{\mathcal{M}}_{min}(\kappa) & = &\overline{{\mathcal{M}}_{min}(\kappa_1)+{%
\mathcal{M}}_{min}(\kappa_2)}^{w^*} \\
&\subseteq & \overline{{\mathcal{M}}_{max}(\kappa_1)+{\mathcal{M}}%
_{max}(\kappa_2)}^{w^*} \subseteq{\mathcal{M}}_{max}(\kappa).
\end{eqnarray*}
Let $T_{\mu}$ be a pseudointegral operator in ${\mathcal{M}}_{min}(\kappa)$
corresponding to a measure $\mu$.
We define measures $\mu_1$ and $\mu_2$ by setting
$$
\mu_1(E)=\mu(E\cap\kappa_1\cap\kappa_2^c), \ \ \mu_2(E)=\mu(E\cap\kappa_2),
\ \ E\subseteq X\times Y \mbox{
Borel}.
$$
It is obvious that $\mu=\mu_1+\mu_2$, hence $T_{\mu} =
T_{\mu_1}+T_{\mu_2}$, and that $\mu_i$ is supported on
$\kappa_i$, $i=1,2$. Thus ${\mathcal{M}}_{min}(\kappa)\subseteq\overline{{%
\mathcal{M}}_{min}(\kappa_1)+{\mathcal{M}}_{min}(\kappa_2)}^{w*}$.

It is clear that ${\mathcal{M}}_{min}(\kappa _{i})\subseteq {\mathcal{M}}%
_{min}(\kappa )$, $i=1,2$, so $\overline{{\mathcal{M}}_{min}(\kappa _{1})+{%
\mathcal{M}}_{min}(\kappa _{2})}^{w\ast }$ $\subseteq $ ${\mathcal{M}}%
_{min}(\kappa )$; thus the first equality is proved.


The second inclusion is obvious. For the last inclusion, note that ${%
\mathcal{M}}_{max}(\kappa _{i})$ $\ \subseteq $ $\ {\mathcal{M}}%
_{max}(\kappa )$, $i=1,2$. Since ${\mathcal{M}}_{max}(\kappa )$ is an
ultraweakly closed space, we conclude that $\overline{{\mathcal{M}}%
_{max}(\kappa _{1})+{\mathcal{M}}_{max}(\kappa _{2})}^{w^{\ast }}\subseteq {%
\mathcal{M}}_{max}(\kappa )$.

Since $\kappa _{1}$ and $\kappa _{2}$ are synthetic, we have that ${\mathcal{%
M}}_{min}(\kappa _{i})={\mathcal{M}}_{max}(\kappa _{i})$, $i=1,2$. If $%
\overline{{\mathcal{M}}_{max}(\kappa _{1})+{\mathcal{M}}_{max}(\kappa _{2})}=%
{\mathcal{M}}_{max}(\kappa )$, then
\begin{eqnarray*}
{\mathcal{M}}_{min}(\kappa ) &=&\overline{{\mathcal{M}}_{min}(\kappa _{1})+{%
\mathcal{M}}_{min}(\kappa _{2})}^{w^{\ast }} \\
&=&\overline{{\mathcal{M}}_{max}(\kappa _{1})+{\mathcal{M}}_{max}(\kappa
_{2})}^{w^{\ast }}={\mathcal{M}}_{max}(\kappa ).
\end{eqnarray*}
Conversely, if $\kappa $ is synthetic, then ${\mathcal{M}}_{min}(\kappa )={%
\mathcal{M}}_{max}(\kappa )$ and thus from the established inclusions it
follows that $\overline{{\mathcal{M}}_{max}(\kappa _{1})+{\mathcal{M}}%
_{max}(\kappa _{2})}={\mathcal{M}}_{max}(\kappa )$. $\diamondsuit $\hfill
\bigskip

We have the following immediate corollaries. The second is a classical
result of Kadison and Singer \cite{ks}.

\begin{corollary}
\label{c_u_all} If $\kappa_1,\kappa_2\subseteq X\times Y$ and $%
\kappa_1\cup\kappa_2= X\times Y$, then ${\mathcal{M}}_{max}(\kappa_1)+{%
\mathcal{M}}_{max}(\kappa_2)$ is dense in ${\mathcal{B}(\mathcal{H}_{1},
\mathcal{H}_{2})}$ in the ultraweak topology.
\end{corollary}

\begin{corollary}
\label{c_u_nest} If ${\mathcal{A}}$ is a nest algebra on a separable Hilbert
space ${\mathcal{H}}$, then ${\mathcal{A}}+{\mathcal{A}}^*$ is dense in ${%
\mathcal{B}}({\mathcal{H}})$ in the ultraweak topology.
\end{corollary}

\bigskip

Since each normalizing masa-bimodule is synthetic \cite{kt}, Proposition \ref
{p_norm_thin} gives a family of synthetic thin masa-bimodules. As the next
example shows, in the class of thin masa-bimodules, there are synthetic
bimodules which are not normalizing.

\begin{example}
\label{e_nn_sy_th} Let ${\mathcal{H}}=L^{2}(0,1)$ (Lebesgue measure) and $%
\kappa =\{(x,y)\in \lbrack 0,1]\times \lbrack 0,1]:y=x\mbox{ or }y=\sqrt{x}%
\} $. Then ${\mathcal{U}}={\mathcal{M}}_{max}(\kappa )$ is a synthetic
non-normalizing thin masa-bimodule.
\end{example}

\noindent \textit{Proof. } Normalizing masa-bimodules have supports of the
form $\{(x,y):f(x)=g(y)\}$ for some Borel functions $f$ and $g$ \cite{kt}.
Thus a support $\lambda $ of a normalizing masa-bimodule has necessarily the
``rectangle'' property
$$
(x,y),(x^{\prime },y),(x,y^{\prime })\in \lambda \Rightarrow (x^{\prime
},y^{\prime })\in \lambda
$$
for marginally almost all points of $\lambda $. (In fact, this property is
characteristic for the supports of normalizing masa-bimodules \cite{t}.) It
is obvious that the set $\kappa $ does not have this property.

To show that ${\mathcal{U}}$ is a thin bimodule, as in Theorem \ref
{e_non-synth}, we use Proposition \ref{p_map_supp}. In this case we consider the
images under the map of ${\mathcal{U}}$ of the sets $\alpha
_{i}^{1}=[0,r_{i}]$ and $\alpha _{i}^{2}=[r_{i},1]$ (see the proof of
Theorem \ref{e_non-synth}).

Now we prove that ${\mathcal{U}}$ is a synthetic bimodule. Let $%
\{a_{i}\}_{i=-\infty }^{+\infty }$ be a sequence such that $a_{i}\in \lbrack
0,1]$, $i\in \mathbb{Z}$, $\{a_{i}\}_{i\leq 0}$ is decreasing and tends to
zero, $\{a_{i}\}_{i\geq 0}$ is increasing and tends to one, $a_{0}=\frac{1}{2%
}$, $a_{i}>a_{i+1}^{2}$, $i\in \mathbb{Z}$. Let $\alpha _{i}=\beta
_{i}^{1}=[a_{i+1},a_{i}]$, $\beta _{i}^{2}=[a_{i+1}^{\frac{1}{2}},a_{i}^{%
\frac{1}{2}}]$, $i\in \mathbb{Z}$. Set $P_{i}=P(\alpha _{i})$ and $%
G_{n}=P((\cup _{i=1}^{n}\alpha _{i})^{c})$, $Q_{i}^{j}=P(\beta _{i}^{j})$, $%
T_{i}^{j}=Q_{i}^{j}TP_{i}$, $j=1,2$. We have that, for each $n\in \mathbb{N}$%
, $T=\sum_{-n\leq i\leq n}TP_{i}+TG_{n}$. Since $T$ is supported on $\kappa $%
, we have that
$$
TP_{i}=(Q_{i}^{1}\vee
Q_{i}^{2})TP_{i}=Q_{i}^{1}TP_{i}+Q_{i}^{2}TP_{i}=T_{i}^{1}+T_{i}^{2}.
$$
Thus, for each $n\in \mathbb{N}$, we have that
$$
T=\sum_{-n\leq i\leq n}(T_{i}^{1}+T_{i}^{2})+TG_{n}.
$$
Now let $S_{n}=\sum_{-n\leq i\leq n}T_{i}^{1}+TG_{n}$. It is easy to see
that the sequence $\{S_{n}\}_{n=1}^{\infty }$ is bounded in norm and thus
there is a subsequence $\{S_{n^{\prime }}\}$ which converges to some
operator $T_{1}$ in the ultraweak topology. Now, the operator $T_{1}$ is
supported on the set
$$
\lambda _{n}=\kappa _{2}\cup \{(x,x):a_{-n}\leq x\leq a_{n}\}.
$$
So $T_{1}$ is supported on the intersection of $\lambda _{n}$, which is
equal to $\kappa _{2}$.

Letting $T_{2}=T-T_{1}$, we have that $\sum_{i=1}^{n}T_{i}^{2}%
\longrightarrow _{n}T_{2}$ in the ultarweak topology and obviously $T_{2}$
is supported on $\kappa _{1}$ (we have that $\sum_{i=1}^{n}T_{i}^{2}$ is
supported on $\kappa _{1}$ for each $n$). Thus we showed that ${\mathcal{U}}=%
{\mathcal{M}}_{max}(\kappa _{1})+{\mathcal{M}}_{max}(\kappa _{2})$. Since
the sets $\kappa _{1}$ and $\kappa _{2}$ are supports of normalizing
masa-bimodules, they are synthetic \cite{kt}. From Proposition \ref{p_u_s}
it follows that ${\mathcal{U}}$ is synthetic. $\diamondsuit $\hfill

Let us note that the minimal normalizing masa-bimodule containing the
bimodule from Example \ref{e_nn_sy_th} has support $\tilde{\kappa}%
=\{(x,y):\exists n\in \mathbb{Z}\mbox{ with }y=x^{2^{n}}\}$.

\bigskip

\textbf{Acknowledgements} A part of the present work was done while the
author was working on his PhD thesis at the University of Athens. The author
wishes to express his deepest gratitude to his supervisor Aristides
Katavolos for his help during the preparation of the paper.
The author would also like to thank the referee for his/her comments.

This work was supported by the Hellenic State Scholarship Foundation.

\textbf{Note} After this work was completed, we were informed that
V. S. Shulman and L. Turowska have independently obtained related
results on the synthesis of finite width masa-bimodules (Section
2); see \cite{st} and \cite{st1}.

2000 Mathematics Subject Classification Numbers:
Primary  47L05; Secondary 47L35

\end{document}